\documentclass{ article}

\begin{document}
\vspace*{.5cm}
\begin{center}
{\Large{\bf   Biharmonic Riemannian maps}}\\
\vspace{.5cm}
 { Bayram \d{S}ahin} \\
\end{center}

\vspace{.5cm}
\begin{center}
{\it Inonu University, Department of Mathematics, 44280,
Malatya-Turkey. E-mail:bsahin@inonu.edu.tr}
\end{center}
\vspace{.5cm}

\noindent {\bf Abstract.} {\small We give necessary and sufficient conditions  for Riemannian maps
to be biharmonic. We also define pseudo umbilical Riemannian maps
as a generalization of pseudo-umbilical submanifolds and show that
such Riemannian maps put some restrictions on the base manifolds.\\

 \noindent{\bf 2000 Mathematics Subject Classification:}
53C43,57R35,53C20.\\

\noindent{\bf Keywords:}Riemannian map, pseudo
umbilical Riemannian map, harmonic map, biharmonic map, Riemannian
submersion.

\pagestyle{myheadings}

\section*{1.Introduction}

  \setcounter{equation}{0}
\renewcommand{\theequation}{1.\arabic{equation}}

Smooth maps between Riemannian manifolds are useful for comparing
geometric structures between two manifolds. Isometric immersions
(Riemannian submanifolds) are basic such maps between Riemannian
manifolds and they are characterized by their Riemannian metrics
and Jacobian matrices. More precisely, a smooth map $F:(M_1,
g_1)\longrightarrow (M_2, g_2)$ between Riemannian manifolds
$(M_1, g_1)$ and $(M_2, g_2)$ is called an isometric immersion if
$F_*$ is injective and
\begin{equation}
g_2(F_*X, F_*Y)=g_1(X, Y)\label{eq:1.1}
\end{equation} for $X, Y $
vector fields tangent to $M_1$, here $F_*$ denotes the derivative
map.

On the other hand, Riemannian submersions between Riemannian manifolds were initiated  by B. O'Neill \cite{O'Neill} and A. Gray
\cite{Gray}, see also \cite{Falcitelli} and \cite{Yano-Kon}. A
smooth map $F:(M_1, g_1)\longrightarrow (M_2, g_2)$ is called
Riemannian submersion if $F_*$ is onto and it satisfies the
equation(\ref{eq:1.1}) for vector fields tangent to the
horizontal space $(kerF_*)^\perp$. For Riemannian submersions between various manifolds, see: \cite{Falcitelli} and
\cite{Yano-Kon}.

In 1992, Fischer  introduced Riemannian maps between Riemannian
manifolds in \cite{Fischer} as a generalization of the notions of
isometric immersions and Riemannian submersions. Let $F:(M_1,
g_1)\longrightarrow (M_2, g_2)$  be a smooth map between
Riemannian manifolds such that $0<rank F<min\{ m, n\}$, where
$dimM_1=m$ and $dimM_2=n$. Then we denote the kernel space of
$F_*$ by $kerF_*$ and consider the orthogonal complementary space
$\mathcal{H}=(ker F_*)^\perp$ to $kerF_*$. Then the tangent bundle
of $M_1$ has the following decomposition

$$TM_1=kerF_* \oplus\mathcal{H}.$$

We denote the range of $F_*$ by $rangeF_*$ and consider the
orthogonal complementary space $(range F_*)^\perp$ to $rangeF_*$
in the tangent bundle $TM_2$ of $M_2$. Since $rankF<min\{ m, n\}$,
we always have $(range F_*)^\perp \neq \{0\}$. Thus the tangent
bundle $TM_2$ of $M_2$ has the following decomposition
$$TM_2=(rangeF_*)\oplus (rangeF_*)^\perp.$$
 Now, a smooth
map $F:(M^{^m}_1,g_1)\longrightarrow (M^{^n}_2, g_2)$ is called a
Riemannian map at $p_1 \in M$ if the horizontal restriction
$F^{^h}_{*p_1}: (ker F_{*p_1})^\perp \longrightarrow (range
F_{*p_1})$  is a linear isometry between the inner product spaces
$((ker F_{*p_1})^\perp, g_1(p_1)\mid_{(ker F_{*p_1})^\perp})$ and
$(range F_{*p_1}, g_2(p_2)\mid_{(range F_{*p_1})})$, $p_2=F(p_1)$.
Therefore Fischer stated in \cite{Fischer} that a Riemannian map
is a map which is as isometric as it can be. In other words, $F_*$
satisfies the equation(\ref{eq:1.1}) for $X, Y$ vector fields
tangent to $\mathcal{H}=(ker F_*)^\perp$. It follows that
isometric immersions and Riemannian submersions are particular
Riemannian maps with $kerF_*=\{ 0 \}$ and $(range F_*)^\perp=\{ 0
\}$. It is known that a Riemannian map is a subimmersion. One of
the main properties of Riemannian maps is that Riemannian maps
satisfy
 the eikonal equation which is a link between geometric optics and physical optics. For Riemannian maps and their applications, see: \cite{Garcia-Kupeli}.

A map between Riemannian manifolds is harmonic if the divergence
of its differential vanishes. Harmonic maps between Riemannian
manifolds provide a rich display of both differential geometric
and analytic phenomena, and they are closely related to the theory
of stochastic processes and to the theory of liquid crystals in
material science.  On the other hand, the biharmonic maps are the
critical points of the bienergy functional and, from this point of
view, generalize harmonic maps. The notion of biharmonic map was
suggested by Eells and Sampson \cite{Eels}. The first variation
formula and, thus, the Euler-Lagrange equation associated to the
bienergy was obtained by Jiang in \cite{Jiang1}, \cite{Jiang2}.
But biharmonic maps have been extensively studied in the last
decade and there are two main research directions. In differential
geometry, many authors have obtained classification results and
constructed many examples. Biharmonicity of immersions was
obtained in \cite{Chen}, \cite{Cadde}, \cite{Oniciuc} and
biharmonic Riemannian submersions were studied in \cite{Oniciuc},
for a survey on biharmonic maps, see:\cite{survey}. From the
analytic point of view, biharmonic maps are solutions of fourth
order strongly elliptic semilinear partial differential equations.
It is known that plane elastic problems can be expressed in terms
of the biharmonic equation. On the other hand, the wave maps are
harmonic maps on Minkowski spaces and the biwave maps are
biharmonic maps on Minkowski spaces.  The wave maps arise in the
analysis of the more difficult hyperbolic Yang-Mills equations
either as special cases or as equations for certain families of
gauge transformations. Such equations arise  in general relativity
for spacetimes with two Killing vector fields. Bi-Yang-Mills
fields, which generalize Yang-Mills fields, have been introduced
by Bejan and Urakawa \cite{Bejan} recently. For relations between
the biwave maps and the bi-Yang-Mills equations, see \cite{Ichima}
and \cite{Jen}. Moreover, in geometric optics\cite{Dietrich}, one
can obtain the eikonal equation by using the wave equation.

In this paper, we mainly investigate the biharmonicity of
Riemannian maps from Riemannian manifolds to space forms. In
section 2, we introduce notations and give fundamental formulas of
the bitension field, then  we obtain some preparatory results of
Riemannian maps in section 3. We also define pseudo umbilical
Riemannian maps as a generalization of pseudo umbilical isometric
immersions, obtain a necessary and sufficient condition for a
Riemannian map to be pseudo umbilical and give a method how to
construct pseudo-umbilical Riemannian maps.  In section 4, we find
necessary and sufficient conditions for Riemannian maps  to be
harmonic and observe that pseudo-umbilical Riemannian maps from
Riemannian manifolds $M_1$ to space forms $M_2(c)$ with additional
conditions must be either harmonic or $c> 0$.

\section*{2.Preliminaries}
  \setcounter{equation}{0}
\renewcommand{\theequation}{2.\arabic{equation}}
In this section we recall some basic materials from \cite{Wood-Baird} and \cite{survey}. Let $(M, g_{_M})$ be a Riemannian manifold and $\mathcal{V}$ be a
$q-$ dimensional distribution on $M.$ Denote its orthogonal
distribution $\mathcal{V}^{\perp}$ by $\mathcal{H}.$ Then, we have
\begin{equation}
TM=\mathcal{V}\oplus \mathcal{H}. \label{eq:2.1}
\end{equation}
$\mathcal{V}$ is called the vertical distribution and
$\mathcal{H}$ is called the horizontal distribution. We use the
same letters to denote the orthogonal projections onto these
distributions.

By the unsymmetrized second fundamental form of $\mathcal{V},$ we
mean the tensor field $A^{\mathcal{V}}$ defined by
\begin{equation}
A^{\mathcal{V}}_E F=\mathcal{H}(\nabla_{\mathcal{V}E} \mathcal{V}
F),\quad  E, F \in \Gamma(TM), \label{eq:2.2}
\end{equation}
where $\nabla$ is the Levi-Civita connection on $M.$ The
symmetrized second fundamental form $B^{\mathcal{V}}$ of
$\mathcal{V}$ is given by
\begin{equation}
B^{\mathcal{V}}(E, F)=\frac{1}{2} \{ A^{\mathcal{V}}_E F+
A^{\mathcal{V}}_F E\}=\frac{1}{2}\{
\mathcal{H}(\nabla_{\mathcal{V}E} \mathcal{V}
F)+\mathcal{H}(\nabla_{\mathcal{V}F} \mathcal{V}
E)\}\label{eq:2.3}
\end{equation}
for  any $ E, F \in \Gamma(TM).$ The integrability tensor of
$\mathcal{V}$ is the tensor field  $I^{\mathcal{V}}$ given by
\begin{equation}
I^{\mathcal{V}}(E, F)=A^{\mathcal{V}}_E F-A^{\mathcal{V}}_F
E-\mathcal{H}([\mathcal{V}E, \mathcal{V}F]). \label{eq:2.4}
\end{equation}
Moreover, the mean curvature vector field of $\mathcal{V}$ is
defined by
\begin{equation}
\mu^{\mathcal{V}}=\frac{1}{q}Trace
B^{\mathcal{V}}=\frac{1}{q}\sum^q_{i=1} \mathcal{H}(\nabla_{e_r}
e_r), \label{eq:2.5}
\end{equation}
where $\{ e_1,..., e_{q}\}$ is a local frame of $\mathcal{V}.$  By
reversing the roles of $\mathcal{V},$ $\mathcal{H},$
$B^{\mathcal{H}},$ $A^{\mathcal{H}}$ and $I^{\mathcal{H}}$ can be
defined similarly. For instance, $B^{\mathcal{H}}$ is defined by
\begin{equation}
B^{\mathcal{H}}(E, F)= \frac{1}{2}\{
\mathcal{V}(\nabla_{\mathcal{H}E} \mathcal{H}
F)+\mathcal{V}(\nabla_{\mathcal{H}F} \mathcal{H} E)\}
\label{eq:2.6}
\end{equation}
and, hence we have
\begin{equation}
\mu^{\mathcal{H}}=\frac{1}{m-q}Trace B^{\mathcal{H}}=\frac{1}{m-q}
\sum^{m-q}_{s=1} \mathcal{V}(\nabla_{E_s} E_s), \label{eq:2.7}
\end{equation}
where $E_1,..., E_{m-q}$ is a local frame of $\mathcal{H}$. A
distribution $\mathcal{D}$ on $M$ is said to be  minimal if, for
each $x \in M$, the mean curvature vector field vanishes.

Let $(M, g_{_M})$ and $(N, g_{_N})$ be Riemannian manifolds and
suppose that $\varphi: M\longrightarrow N$ is a smooth map between
them. Then the differential ${\varphi}_*$ of $\varphi$ can be
viewed a section of the bundle $Hom(TM,
\varphi^{-1}TN)\longrightarrow M,$ where $\varphi^{-1}TN$ is the
pullback bundle which has fibres
$(\varphi^{-1}TN)_p=T_{\varphi(p)} N, p \in M.$ $Hom(TM,
\varphi^{-1}TN)$ has a connection $\nabla$ induced from the
Levi-Civita connection $\nabla^M$ and the pullback connection.
Then the second fundamental form of $\varphi$ is given by
\begin{equation}
(\nabla {\varphi}_*)(X, Y)=\nabla^{\varphi}_X
{\varphi}_*(Y)-{\varphi}_*(\nabla^M_X Y) \label{eq:2.10}
\end{equation}
for $X, Y \in \Gamma(TM).$ It is known that the second fundamental
form is symmetric. Let $\varphi: (M, g_{_M}) \longrightarrow (N,
g_{_N})$ be a smooth map between Riemannian manifolds and assume
$M$ is compact, then its energy is
\begin{eqnarray}
E(\varphi)=\int_{M}e(\varphi)\,v_{g}=\frac{1}{2}\int_{M}|d\varphi|^{2}\,v_{g}.\nonumber
\end{eqnarray}
The critical points of $E$  are called harmonic maps. Standard
arguments yield the associated Euler-Lagrange equation, the
vanishing of the tension field:$\tau(\varphi)=trace (\nabla
{\varphi}_*)$. Let $\varphi: (M, g_{_M}) \longrightarrow (N,
g_{_N})$ be a smooth map between Riemannian manifolds. Define its
bienergy as

$$E^2(\varphi)=\frac{1}{2}\int_{M}|\tau(\varphi)|^{2}\,v_{g}.$$
Critical points of the functional $E^2$ are called biharmonic maps
and its associated Euler-Lagrange equation is the vanishing of the
bitension field
\begin{equation}
\tau^2(\varphi)=-\Delta^{\varphi}\tau(\varphi)-trace_{g_{_M}}R^N(d\varphi,\tau(\varphi))d\varphi,
\label{eq:2.11}
\end{equation}
where
$\Delta^{\varphi}\tau(\varphi)=-trace_{g_{_M}}(\nabla^{\varphi}\nabla^{\varphi}-\nabla^{\varphi}_{_{\nabla}})$
is the Laplacian on the sections of $\varphi^{-1}(TN)$ and $R^N$
is the Riemann curvature operator on $(N, g_{_N})$. A map between
two Riemannian manifolds is said to be proper biharmonic if it is
a non-harmonic biharmonic map.

\section*{3.Riemannian maps}
  \setcounter{equation}{0}
\renewcommand{\theequation}{3.\arabic{equation}}

In this section, we obtain some new results which will be using in
the next section. First note that in \cite{Sahin2} we showed that
the second fundamental form  $(\nabla F_*)(X,Y)$, $\forall X, Y
\in \Gamma((ker F_*)^\perp)$, of a Riemannian map has no
components in $range F_*$.

\noindent{\bf Lemma~3.1.~}{\it
Let $F$ be a Riemannian map from a Riemannian manifold $(M_1,g_1)$
to a Riemannian manifold $(M_2,g_2)$. Then}
\begin{equation}
g_2((\nabla F_*)(X,Y),F_*(Z))=0, \forall X,Y,Z \in \Gamma((ker
F_*)^\perp). \label{eq:3.1}
\end{equation}

As a result of Lemma 3.1, we have
\begin{equation}
(\nabla F_*)(X,Y) \in \Gamma((range F_*)^\perp), \forall X,Y \in
\Gamma((ker F_*)^\perp). \label{eq:3.4}
\end{equation}

Also from \cite{Sahin1}, we have the following.

\noindent{\bf Lemma~3.2~}{\it Let $F: (M, g_{_M})
\longrightarrow (N,g_{_N})$ be a Riemannian map between Riemannian
manifolds. Then the tension field $\tau$ of $F$ is}
\begin{equation}
\tau=-m_1F_*(\mu^{ker F_*})+m_2H_2, \label{eq:3.5}
\end{equation}
{\it where $m_1=dim(ker F_*), m_2=rank F$, $\mu^{ker F_*}$ and $H_2$
are the mean curvature vector fields of the distributions of $ker
F_*$ and $range F_*$, respectively.}

From now on, for simplicity, we denote by $\nabla^2$ both the
Levi-Civita connection of $(M_2, g_2)$ and its pullback along $F$.
 Then according to \cite{Nore}, for any vector field $X$ on $M_1$ and any section $V$ of $(range F_*)^\perp$, where $(range F_*)^\perp$ is the
 subbundle of $F^{-1}(TM_2)$ with fiber $(F_*(T_pM))^\perp$-orthogonal complement of $F_*(T_pM)$ for $g_2$ over $p$, we have
 $\nabla^{^{F \perp}}_XV$ which is the orthogonal projection of $\nabla^2_XV$ on $(F_*(TM))^\perp$. In \cite{Nore}, the author also showed that $\nabla^{^{F \perp}}$ is a linear connection on $(F_*(TM))^\perp$ such that $\nabla^{^{F \perp}}g_2=0$. We now  define $A_{V}$ as
\begin{equation}
\nabla^2_{_{F_*X}}V=-A_{_V}F_*X+\nabla^{^{F \perp}}_{_{X}}V,
\label{eq:3.6}
\end{equation}
where $A_{_V}F_*X$ is the tangential component (a vector field
along $F$) of $\nabla^2_{_{F_*X}}V$. It is easy to see that $A_V
F_*X$ is bilinear in $V$ and $F_*X$ and $A_V F_*X$ at $p$ depends
only on $V_p$ and $F_{*p}X_p$. By direct computations, we obtain
\begin{equation}
g_2(A_{_V} F_*X,F_*Y)=g_2(V, (\nabla F_*)(X,Y)), \label{eq:3.7}
\end{equation}
for $X, Y \in \Gamma((ker F_*)^\perp)$ and $V \in \Gamma((range
F_*)^\perp)$. Since $(\nabla F_*)$ is symmetric, it follows that
$A_{_V}$ is a symmetric linear transformation of $range F_*$.

We now define pseudo-umbilical Riemannian maps as a generalization
of pseudo-umbilical isometric immersions. Pseudo-umbilical
Riemannian maps will be useful when we deal with the biharmonicity
of Riemannian maps.

\noindent{\bf Definitioñ~3.1.~} Let $F:(M,g_1) \longrightarrow (M_2,g_2)$ be a Riemannian
map between Riemannian manifolds $M_1$ and $M_2$. Then we say
that $F$ is a pseudo-umbilical Riemannian map if
\begin{equation}
A_{H_2}F_*(X)=\lambda F_*(X) \label{eq:3.8}
\end{equation}
for $\lambda \in C^{\infty}(M_1)$ and $X \in \Gamma((ker
F_*)^\perp)$.

Here we present an useful formula for pseudo umbilical Riemannian
maps by using(\ref{eq:3.7}) and(\ref{eq:3.8}).

\noindent{\bf Proposition~3.1.~}{\it Let $F:(M,g_1) \longrightarrow (M_2,g_2)$ be a Riemannian map
between Riemannian manifolds $M_1$ and $M_2$. Then $F$ is
pseudo-umbilical if and only if}
\begin{equation}
g_2((\nabla F_*)(X,Y),H_2)=g_1(X,Y)g_2(H_2,H_2) \label{eq:3.9}
\end{equation}
 {\it for $X, Y \in \Gamma((ker F_*)^\perp)$.}

\noindent{\bf Proof.~} Let $\{\tilde{e}_1,...,\tilde{e}_{m_1},e_{1},...,e_{m_2}\}$
be an orthonormal basis of $\Gamma(TM_1)$ such that
$\{\tilde{e}_1,...,\tilde{e}_{m_1}\}$ is an orthonormal basis of
$ker F_*$ and $\{e_{1},...,e_{m_2}\}$ is an orthonormal basis of
$(ker F_*)^\perp$. Then  since $F$ is a Riemannian map we have
$$\sum^{m_2}_{i=1}g_2(A_{H_2}F_*(e_i),F_*(e_i))=m_2\lambda.$$
Using(\ref{eq:3.7}), we get
$$\sum^{m_2}_{i=1}g_2(\frac{1}{m_2}(\nabla F_*)(e_i,e_i),H_2)=\lambda.$$
Thus we obtain
\begin{equation}
\lambda=g_2(H_2,H_2). \label{eq:3.10}
\end{equation}
Then, from (\ref{eq:3.7}), (\ref{eq:3.8}) and(\ref{eq:3.10}) we
have(\ref{eq:3.9}). The converse is clear.\\

It is known that the composition of a Riemannian submersion and an
isometric immersion is a Riemannian map \cite{Fischer}.Using this
we have the following.

\noindent{\bf Theorem~3.1.~}{\it Let $F_1:(M_1,g_1)\longrightarrow (M_2,g_2)$ be a Riemannian
submersion and $F_2:(M_2,g_2)\longrightarrow (M_3,g_3)$ a
pseudo-umbilical isometric immersion. Then the map $F_2 \circ F_1$
is a pseudo umbilical Riemannian map.}\\

\noindent{\bf Proof.~} From the second fundamental form of the composite map $F_2
\circ F_1$ \cite{Wood-Baird}, we have
$$(\nabla (F_2 \circ F_1)_*(X,Y)=F_{2*}((\nabla
F_{1*})(X,Y))+(\nabla F_{2*})(F_{1*}X,F_{1*}Y)$$ for $X, Y \in
\Gamma((ker F_{1*})^\perp)$. Then proof follows from the
definition of pseudo-umbilical submanifolds.\\

\noindent{\bf Remark.~3.1.~} We note that above theorem gives a method to find examples of
pseudo umbilical Riemannian maps. It also tells that if one has an
example of pseudo-umbilical submanifolds, it is possible to find
an example of pseudo umbilical Riemannian maps. For examples of
pseudo umbilical submanifolds, see: \cite{Chenpseudo}.

\section*{4.Biharmonicity of Riemannian maps}
  \setcounter{equation}{0}
\renewcommand{\theequation}{4.\arabic{equation}}

In this section we obtain the biharmonicity of Riemannian maps
between Riemannian manifolds. We also show that pseudo-umbilical
biharmonic Riemannian maps put some restrictions on the total
manifold of such maps.

Let $F:(M_1,g_1) \longrightarrow (M_2,g_2)$ be a map between
Riemannian manifolds $(M_1,g_1)$ and $(M_2,g_2)$.  Then the
adjoint map ${^*F}_*$ of $F_*$ is characterized by
$g_1(x,{^*F}_{*p_1}y)=g_2(F_{*p_1}x,y)$ for $x \in T_{p_1}M_1$, $y
\in T_{F(p_1)}M_2$ and $p_1 \in M_1$. Considering $F^{h}_*$ at
each $p_1 \in M_1$ as a linear transformation
$$F^{h}_{*{p_1}}: ((ker F_*)^{\perp}(p_1),{g_1}_{{p_1}((ker
F_*)^\perp (p_1))}) \longrightarrow (range
F_*(p_2),{g_2}_{{p_2}(range F_*)(p_2))}),$$ we will denote the
adjoint of $F^{h}_* $ by ${^*F^{h}}_{*p_1}.$ Let ${^*F}_{*p_1}$ be
the adjoint of $F_{*p_1}:(T_{p_1}M_1,{g_1}_{p_2}) \longrightarrow
(T_{p_2}M_2, {g_2}_{p_2})$. Then the linear transformation
$$({^*F}_{*p_1})^h:range F_*(p_2) \longrightarrow (ker F_*)^{\perp}(p_1)$$
defined by $({^*F}_{*p_1})^hy= {^*F}_{*p_1}y$, where $y \in
\Gamma(range F_{*p_1}), p_2=F(p_1)$, is an isomorphism and
$(F^h_{*p_1})^{-1}=({^*F}_{*p_1})^h={^*(F^h_{*p_1})}$.

We also recall that the curvature tensor $R$ of a space form
$(M(c),g)$ is given by
\begin{equation}
R(X,Y)Z=c\{g(Y,Z)X-g(X,Z)Y\}. \label{eq:4.1}
\end{equation}
 We are now ready to prove the following theorem which gives necessary
 and sufficient conditions for a Riemannian map to be biharmonic.

\noindent{\bf Theorem~~4.1.~} {\it Let $F$ be a Riemannian map from a Riemannian manifold
$(M_1,g_1)$ to a space form $(M_2(c),g_2)$. Then $F$ is biharmonic
if and only if}
\begin{eqnarray}
&&m_1trace A_{(\nabla F_*)(.,\mu^{ker F_*})}F_*(.)-m_1trace
F_*(\nabla_{(.)} \nabla_{(.)} \mu^{ker F_*})\nonumber\\
&&-m_2trace F_*(\nabla_{(.)}{^*F}_*(A_{H_2}F_*(.)))- m_2trace
A_{\nabla^{^{F \perp}}_{F_*(.)}H_2}F_*(.)\nonumber\\
&&-m_1c(m_2-1)F_*(\mu^{ker F_*})=0 \label{eq:4.2}
\end{eqnarray}
{\it and}
\begin{eqnarray}
&&m_1trace \nabla^{^{F \perp}}_{F_*(.)}(\nabla F_*)(.,\mu^{ker
F_*})+m_1trace(\nabla F_*)(.,\nabla_{(.)}\mu^{ker F_*})\nonumber\\
&&+m_2trace (\nabla
F_*)(.,{^*F}_*(A_{H_2}F_*(.)))-m_2\Delta^{R^\perp}H_2\nonumber\\
&&-m^2_2cH_2=0.\label{eq:4.3}
\end{eqnarray}

\noindent{\bf Proof.~}First of all, from(\ref{eq:4.1}) and
(\ref{eq:3.6}) we have
\begin{equation}
trace R^2(F_*(.),\tau(F))F_*(.)=m_1c(m_2-1)F_*(\mu^{ker
F_*})-m_2^2cH_2, \label{eq:4.4}
\end{equation}
where $R^2$ is the curvature tensor field of $M_2$. Let
$\{\tilde{e}_1,...,\tilde{e}_{m_1},e_{1},...,e_{m_2}\}$ be a local
orthonormal frame on $M_1$, geodesic at $p \in M_1$ such that
$\{\tilde{e}_1,...,\tilde{e}_{m_1}\}$ is an orthonormal basis of
$ker F_*$ and $\{e_{1},...,e_{m_2}\}$ is an orthonormal basis of
$(ker F_*)^\perp$. At $p$ we have

\begin{eqnarray}
\Delta
\tau(F)&=&-\sum^{m_2}_{i=1}\nabla^{^F}_{e_i}\nabla^{^F}_{e_i}\tau(F)\nonumber\\
&=&-\sum^{m_2}_{i=1}\nabla^{^F}_{e_i}\{\nabla^{^F}_{e_i}(-m_1F_*(\mu^{ker
F_*})+m_2H_2)\}.\nonumber
\end{eqnarray}
Then using(\ref{eq:2.10}),(\ref{eq:3.4}) and(\ref{eq:3.6}) we
get
\begin{eqnarray}
\Delta \tau(F)&=&-\sum^{m_2}_{i=1}\nabla^{^F}_{e_i}\{-m_1(\nabla
F_*)(e_i,\mu^{ker F_*})-m_1F_*(\nabla_{e_i}\mu^{ker F_*})\nonumber\\
&+&m_2(-A_{H_2}F_*(e_i)+\nabla^{^{F
\perp}}_{F_*(e_i)}H_2)\}.\nonumber
\end{eqnarray}
Using again (\ref{eq:2.10}),(\ref{eq:3.4}) and(\ref{eq:3.6}) we
obtain
\begin{eqnarray}
\Delta \tau(F)&=&m_1\sum^{m_2}_{i=1}-A_{(\nabla
F_*)(e_i,\mu^{kerF_*})}F_*(e_i)+\nabla^{^{F
\perp}}_{F_*(e_i)}(\nabla
F_*)(e_i,\mu^{ker F_*})\nonumber\\
&+&m_1\sum^{m_2}_{i=1}(\nabla F_*)(e_i,\nabla_{e_i}\mu^{ker
F_*})+F_*(\nabla_{e_i}\nabla_{e_i}\mu^{ker F_*})\nonumber\\
&+&m_2\sum^{m_2}_{i=1}\nabla^{^F}_{e_i}A_{H_2}F_*(e_i)-m_2\sum^{m_2}_{i=1}-A_{\nabla^{^{F \perp}}_{F_*(e_i)}H_2}F_*(e_i)\nonumber\\
&+&\nabla^{^{F \perp}}_{F_*(e_i)}\nabla^{^{F
\perp}}_{F_*(e_i)}H_2.\nonumber
\end{eqnarray}
On the other hand, since $A_{H_2}F_*(e_i) \in \Gamma(F_*((ker
F_*)^\perp))$, we can write
$$F_*(X)=A_{H_2}F_*(e_i)$$
for $X \in \Gamma((ker F_*)^\perp)$, where
$$X=(F_*)^{-1}(A_{H_2}F_*(e_i))={^*F}_*(A_{H_2}F_*(e_i)).$$
Then using(\ref{eq:2.10}) we have
$$ \nabla^{^F}_{e_i}A_{H_2}F_*(e_i)=(\nabla
F_*)(e_i,{^*F}_*(A_{H_2}F_*(e_i)))+F_*(\nabla_{e_i}{^*F}_*(A_{H_2}F_*(e_i))).$$
Thus we obtain
\begin{eqnarray}
\Delta \tau(F)&=&m_1\sum^{m_2}_{i=1}-A_{(\nabla
F_*)(e_i,\mu^{kerF_*})}F_*(e_i)+\nabla^{^{F
\perp}}_{F_*(e_i)}(\nabla
F_*)(e_i,\mu^{ker F_*})\nonumber\\
&+&m_1\sum^{m_2}_{i=1}(\nabla F_*)(e_i,\nabla_{e_i}\mu^{ker
F_*})+F_*(\nabla_{e_i}\nabla_{e_i}\mu^{ker F_*})\nonumber\\
&+&m_2\sum^{m_2}_{i=1}(\nabla
F_*)(e_i,{^*F}_*(A_{H_2}F_*(e_i)))+F_*(\nabla_{e_i}{^*F}_*(A_{H_2}F_*(e_i)))\nonumber\\
&-&m_2\sum^{m_2}_{i=1}-A_{\nabla^{^{F
\perp}}_{F_*(e_i)}H_2}F_*(e_i)+\nabla^{^{F
\perp}}_{F_*(e_i)}\nabla^{^{F \perp}}_{F_*(e_i)}H_2.\label{eq:4.5}
\end{eqnarray}
Thus putting (\ref{eq:4.4}) and (\ref{eq:4.5}) in(\ref{eq:2.11})
and then  taking the $F_*((ker F_*)^\perp)=range F_*$ and $(range
F_*)^\perp$ parts we have(\ref{eq:4.2}) and(\ref{eq:4.3}).\\

In particular, we have the following.\\

\noindent{\bf Corollary~4.1.~}{\it Let $F$ be a Riemannian
map from a Riemannian manifold $(M_1,g_1)$ to a space form
$(M_2(c),g_2)$. If the mean curvature vector fields of $range F_*$
and $ker F_*$ are parallel, then $F$ is biharmonic if and only if}
\begin{eqnarray}
&&m_1trace A_{(\nabla F_*)(.,\mu^{ker F_*})}F_*(.)-m_2trace
F_*(\nabla_{(.)}{^*F}_*(A_{H_2}F_*(.))\nonumber\\
&&-m_1c(m_2-1)F_*(\mu^{ker F_*})=0 \nonumber
\end{eqnarray}
{\it and}
\begin{eqnarray}
&&m_1trace \nabla^{^{F \perp}}_{F_*(.)}(\nabla F_*)(.,\mu^{ker
F_*})+m_2trace (\nabla
F_*)(.,{^*F}_*(A_{H_2}F_*(.)))\nonumber\\
&&-m^2_2cH_2=0\nonumber
\end{eqnarray}

We also have the following result for pseudo-umbilical Riemannian
maps.\\

\noindent{\bf Theorem~4.2.~}{\it Let $F$ be a pseudo-umbilical biharmonic Riemannian map
from a Riemannian manifold $(M_1,g_1)$ to a space form
$(M_2(c),g_2)$ such that the distribution $ker F_*$ is minimal and
the mean curvature vector  field  $H_2$ is parallel. Then either
$F$ is harmonic or $c=\parallel H_2
\parallel^2$.}

\noindent{\bf Proof.~} First note that it is easy to see that $\parallel H_2
\parallel^2$ is constant. If $F$ is biharmonic Riemannian map such
that $\mu^{ker F_*}=0$ and $H_2$ is parallel, then from
(\ref{eq:4.3}) we have
$$
m_2\sum^{m_2}_{i=1}(\nabla
F_*)(e_i,{^*F}_*(A_{H_2}F_*(e_i)))-m_2^2cH_2=0.$$ Since $F$ is
pseudo umbilical, we get
$$
m_2\sum^{m_2}_{i=1}(\nabla F_*)(e_i,{^*F}_*(\parallel H_2
\parallel^2F_*(e_i)))-m_2^2cH_2=0.$$
On the other hand, from the linear map ${^*F}_*$ and ${^*F}_*\circ
F_*=I$ (identity map), we obtain
$$
m_2\sum^{m_2}_{i=1}(\nabla F_*)(e_i,\parallel H_2
\parallel^2e_i))-m_2^2cH_2=0.$$
Since the second fundamental form is also linear in its arguments,
it follows that
$$
m_2\parallel H_2
\parallel^2\sum^{m_2}_{i=1}(\nabla F_*)(e_i,e_i))-m_2^2cH_2=0.$$
Hence we have
$$
m^2_2\parallel H_2
\parallel^2H_2-m_2^2cH_2=0$$
which implies that
\begin{equation}
(\parallel H_2
\parallel^2-c)H_2=0.\label{eq:4.6}
\end{equation}
Thus either $H_2=0$ or $(\parallel H_2
\parallel^2-c)=0$. If $H_2=0$, then Lemma 3.2 implies that $F$ is
harmonic, thus proof is complete.\\

From(\ref{eq:4.6}), we have the following result which puts some
restrictions on $M_2(c)$.\\

\noindent{\bf Corollary~4.2.~}{\it There exist no proper biharmonic pseudo umbilical Riemannian
maps from a Riemannian manifold to space forms $(M_2(c)$ with
$c\leq 0$ such that the distribution $ker F_*$ is minimal and the
mean curvature vector field $H_2$ is parallel.}\\

\noindent{\bf Remark~4.1.~}In this paper, we
investigate the biharmonicity of Riemannian maps between
Riemannian manifolds. Our results give some clues to investigate
the biharmonicity of arbitrary maps between Riemannian manifolds.
They also give a method to investigate the geometry of Riemannian
maps. Since Riemannian maps are  solutions of the eikonal
equations which can be obtained starting from the wave equation,
biharmonic maps are solutions of fourth order strongly elliptic
semilinear partial differential equations and they are related to
the biwave equation and bi-Yang-Mills fields, biharmonic
Riemannian maps have potential for further research in terms of
partial differential equations, geometric optics and mathematical
physics.

\end{document}